\documentclass[a4paper,oneside,english]{amsart}
\usepackage[active]{srcltx}
\usepackage{amsthm}
\usepackage{amstext}
\usepackage{amssymb}

\makeatletter

\newcommand{\noun}[1]{\textsc{#1}}

\theoremstyle{plain}
\newtheorem{thm}{\protect\theoremname}

\makeatother

\usepackage{babel}
\providecommand{\theoremname}{Theorem}

\begin{document}

\title[Symmetric properties of the multiple $q$-Euler zeta functions]{A Note on Symmetric properties of the multiple $q$-Euler zeta functions
and higher-order $q$-Euler polynomials}

\author{Dae San Kim \and Taekyun Kim}
\begin{abstract}
Recently, the higher-order $q$-Euler polynomials and multiple $q$-Euler
zeta functions are introduced by T. Kim (\cite{key-8,key-7}). In
this paper, we investigate some symmetric properties of the multiple
$q$-Euler zeta function and derive various identities concerning
the higher-order $q$-Euler polynomials from the symmetric properties
of the multiple $q$-Euler zeta functions.
\end{abstract}
\maketitle

\section{Introduction}

$\,$

For $q\in\mathbb{C}$ with $\left|q\right|<1$, the $q$-number
is defined by $\left[x\right]_{q}=\frac{1-q^{x}}{1-q}$ . Note that
${\displaystyle \lim_{q\rightarrow1}\left[x\right]_{q}}=x$. As is
well known, the Euler polynomials of order $r$$\left(\in\mathbb{N}\right)$
are defined by the generating function to be
\begin{equation}
\left(\frac{2}{e^{t}+1}\right)^{r}e^{xt}=\underset{r\textrm{-times}}{\underbrace{\left(\frac{2}{e^{t}+1}\right)\times\cdots\times\left(\frac{2}{e^{t}+1}\right)}}e^{xt}=\sum_{n=0}^{\infty}E_{n}^{\left(r\right)}\left(x\right)\frac{t^{n}}{n!}.\label{eq:1}
\end{equation}

When $x=0$, $E_{n}^{\left(r\right)}=E_{n}^{\left(r\right)}\left(0\right)$
are called the Euler numbers of order $r$ (see {[}1-13{]}).

In \cite{key-7}, T. Kim considered the $q$-extension of higher-order
Euler polynomials which are given by the generating function to be
\begin{align}
F_{q}^{\left(r\right)}\left(t,x\right) & =\left[2\right]_{q}^{r}\sum_{m_{1},\cdots,m_{r}=0}^{\infty}\left(-q\right)^{m_{1}+\cdots+m_{r}}e^{\left[m_{1}+\cdots+m_{r}+x\right]_{q}t}\label{eq:2}\\
 & =\sum_{n=0}^{\infty}E_{n,q}^{\left(r\right)}\left(x\right)\frac{t^{n}}{n!}.\nonumber
\end{align}

Note that ${\displaystyle \lim_{q\rightarrow1}}F_{q}^{\left(r\right)}\left(t,x\right)=\left(\frac{2}{e^{t}+1}\right)^{r}e^{xt}={\displaystyle \sum_{n=0}^{\infty}E_{n}^{\left(r\right)}\left(x\right)\frac{t^{n}}{n!}}$.

When $x=0$, $E_{n,q}^{\left(r\right)}=E_{n,q}^{\left(r\right)}\left(0\right)$
are called the $q$-Euler numbers of order $r$$\left(\in\mathbb{N}\right)$.

In \cite{key-11}, Rim et al. have studied the properties of $q$-Euler
polynomials due to T. Kim.

From (\ref{eq:2}), we note that
\begin{align}
E_{n,q}^{\left(r\right)}\left(x\right) & =\sum_{l=0}^{n}\dbinom{n}{l}q^{lx}E_{l,q}^{\left(r\right)}\left[x\right]_{q}^{n-l}\label{eq:3}\\
 & =\left(q^{x}E_{q}^{\left(r\right)}+\left[x\right]_{q}\right)^{n},\nonumber
\end{align}
with the usual convention about replacing $\left(E_{q}^{\left(r\right)}\right)^{n}$
by $E_{n,q}^{\left(r\right)}$.

In \cite{key-7}, T. Kim considered the multiple $q$-Euler zeta function
which interpolates higher-order $q$-Euler polynomials at negative
integers as follows :
\begin{eqnarray}
\zeta_{q,r}\left(s,x\right) & = & \left[2\right]_{q}^{r}\sum_{m_{1},\cdots,m_{r}=0}^{\infty}\frac{\left(-q\right)^{m_{1}+\cdots+m_{r}}}{\left[m_{1}+\cdots+m_{r}+x\right]_{q}^{s}}\label{eq:4}\\
 & = & \left[2\right]_{q}^{r}\sum_{m=0}^{\infty}\dbinom{m+r-1}{m}_{q}\left(-q\right)^{m}\frac{1}{\left[m+x\right]_{q}^{s}},\nonumber
\end{eqnarray}
where $s\in\mathbb{C}$ and $x\in\mathbb{R}$ with $x\ne0,\,-1,\,-2,\,-3,\,\cdots$.

By using Cauchy residue theorem and Laurent series, we note that
\begin{equation}
\zeta_{q,r}\left(-n,x\right)=E_{n,q}^{\left(r\right)}\left(x\right),\quad\textrm{where }r\in\mathbb{Z}_{\ge0}.\label{eq:5}
\end{equation}

Recently, D.S. Kim et al. (\cite{key-5}) introduced some interesting
and important symmetric identities of the $q$-Euler polynomials which
are derived from the symmetric properties of $q$-Euler zeta function.
Indeed, their identities are a part of an answer to an open question
for the symmetric identities of Carlitz's type $q$-Euler polynomials
in \cite{key-6}.

In order to find a generalization of identities of D. S. Kim et al.
(\cite{key-5}), we consider symmetric properties of the multiple
$q$-Euler zeta function. From the symmetric properties of multiple
$q$-Euler zeta function, we derive identities of symmetry for the
higher-order $q$-Euler polynomials.

\section{Some identities of higher-order $q$-Euler polynomials}

$\,$

For $a,\, b\in\mathbb{N}$ with $a\equiv1$ (mod $2$) and $b\equiv1$
(mod 2), we observe that
\begin{align}
 & \frac{1}{\left[2\right]_{q^{a}}^{r}}\zeta_{q^{a},r}\left(s,\, bx+\frac{b}{a}\left(j_{1}+\cdots+j_{r}\right)\right)\label{eq:6}\\
= & \sum_{n_{1},\cdots,\, n_{r}=0}^{\infty}\frac{\left(-1\right)^{n_{1}+\cdots+n_{r}}q^{a\left(n_{1}+\cdots+n_{r}\right)}}{\left[n_{1}+\cdots+n_{r}+bx+\frac{b}{a}\left(j_{1}+\cdots+j_{r}\right)\right]_{q^{a}}^{s}}\nonumber \\
= & \left[a\right]_{q}^{s}\sum_{n_{1},\cdots,\, n_{r}=0}^{\infty}\frac{\left(-1\right)^{n_{1}+\cdots+n_{r}}q^{a\left(n_{1}+\cdots+n_{r}\right)}}{\left[a\left(n_{1}+\cdots+n_{r}\right)+abx+b\left(j_{1}+\cdots+j_{r}\right)\right]_{q}^{s}}\nonumber \\
= & \left[a\right]_{q}^{s}\sum_{n_{1},\cdots,\, n_{r}=0}^{\infty}\sum_{i_{1},\cdots,i_{r}=0}^{b-1}\frac{\left(-1\right)^{{\displaystyle {\textstyle \sum_{l=1}^{r}\left(i_{l}+bn_{l}\right)}}}q^{a{\displaystyle {\textstyle \sum_{l=1}^{r}\left(i_{l}+bn_{l}\right)}}}}{\left[ab{\displaystyle {\textstyle \sum_{l=1}^{r}\left(x+n_{l}\right)}+b{\textstyle \sum_{l=1}^{r}j_{l}}+a{\textstyle \sum_{l=1}^{r}i_{l}}}\right]_{q}^{s}}.\nonumber
\end{align}

From (\ref{eq:6}), we note that
\begin{align}
 & \frac{\left[b\right]_{q}^{s}}{\left[2\right]_{q^{a}}^{r}}\sum_{j_{1},\cdots,\, j_{r}=0}^{a-1}\left(-1\right)^{{\textstyle {\displaystyle {\textstyle \sum_{l=1}^{r}j_{l}}}}}q^{b{\displaystyle {\textstyle \sum_{l=1}^{r}j_{l}}}}\zeta_{q^a,r}\left(s,\, bx+\frac{b}{a}\left(j_{1}+\cdots+j_{r}\right)\right)\label{eq:7}\\
= & \left[a\right]_{q}^{s}\left[b\right]_{q}^{s}\sum_{j_{1},\cdots,\, j_{r}=0}^{a-1}\sum_{i_{1},\cdots,\, i_{r}=0}^{b-1}\nonumber \\
 & \times\sum_{n_{1},\cdots,\, n_{r}=0}^{\infty}\frac{\left(-1\right)^{{\displaystyle {\textstyle \sum_{l=1}^{r}\left(i_{l}+j_{l}+n_{l}\right)}}}q^{{\displaystyle {\textstyle \sum_{l=1}^{r}\left(bj_{l}+ai_{l}+abn_{l}\right)}}}}{\left[ab{\displaystyle {\textstyle \sum_{l=1}^{r}\left(x+n_{l}\right)+b{\displaystyle {\textstyle \sum_{l=1}^{r}j_{l}}+a{\displaystyle {\textstyle \sum_{l=1}^{r}i_{l}}}}}}\right]_{q}^{s}}.\nonumber
\end{align}

By the same method as (\ref{eq:7}), we get

\begin{align}
 & \frac{\left[a\right]_{q}^{s}}{\left[2\right]_{q^{b}}^{r}}\sum_{j_{1},\cdots,\, j_{r}=0}^{b-1}\left(-1\right)^{{\displaystyle {\textstyle \sum_{l=1}^{r}j_{l}}}}q^{a{\displaystyle {\textstyle \sum_{l=1}^{r}j_{l}}}}\zeta\left(s,\, ax+\frac{a}{b}\left(j_{1}+\cdots+j_{r}\right)\right)\label{eq:8}\\
= & \left[a\right]_{q}^{s}\left[b\right]_{q}^{s}\sum_{j_{1},\cdots,\, j_{r}=0}^{b-1}\sum_{i_{1},\cdots,\, i_{r}=0}^{a-1}\nonumber \\
 & \times\sum_{n_{1},\cdots,\, n_{r}=0}^{\infty}\frac{\left(-1\right)^{{\displaystyle {\textstyle \sum_{l=1}^{r}\left(i_{l}+j_{l}+n_{l}\right)}}}q^{{\displaystyle {\textstyle \sum_{l=1}^{r}\left(bi_{l}+aj_{l}+abn_{l}\right)}}}}{\left[ab{\displaystyle {\textstyle \sum_{l=1}^{r}\left(x+n_{l}\right)}+a{\textstyle \sum_{l=1}^{r}j_{l}}+b{\displaystyle {\textstyle \sum_{l=1}^{r}i_{l}}}}\right]_{q}^{s}}.\nonumber
\end{align}

Threfore, by (\ref{eq:7}) and (\ref{eq:8}), we obtain the following
theorem.
\begin{thm}
\label{thm:1} For $a,\, b\in\mathbb{N}$ with $a\equiv1$ $\textnormal{(mod }2\mathnormal{)}$
and $b\equiv1$ $\textnormal{(mod }2\mathnormal{)}$, we have
\begin{align*}
 & \left[2\right]_{q^{b}}^{r}\left[b\right]_{q}^{s}\sum_{j_{1},\cdots,j_{r}=0}^{a-1}\left(-1\right)^{{\displaystyle {\textstyle \sum_{l=1}^{r}j_{l}}}}q^{b{\displaystyle {\textstyle \sum_{l=1}^{r}j_{l}}}}\zeta_{q^{a},r}\left(s,bx+\frac{b}{a}\left(j_{1}+\cdots+j_{r}\right)\right)\\
= & \left[2\right]_{q^{a}}^{r}\left[a\right]_{q}^{s}\sum_{j_{1},\cdots,j_{r}=0}^{b-1}\left(-1\right)^{{\displaystyle {\textstyle \sum_{l=1}^{r}j_{l}}}}q^{a{\displaystyle {\textstyle \sum_{l=1}^{r}j_{l}}}}\zeta_{q^{b},r}\left(s,ax+\frac{a}{b}\left(j_{1}+\cdots+j_{r}\right)\right).
\end{align*}

\end{thm}
$\,$

From (\ref{eq:5}) and Theorem \ref{thm:1}, we obtain the following
theorem.
\begin{thm}
For $n\ge0$ and $a,\, b\in\mathbb{N}$ with $a\equiv1$ $\textnormal{(mod }2\mathnormal{)}$
and $b\equiv1$ $\textnormal{(mod }2\mathnormal{)}$, we have

\begin{align*}
 & \left[2\right]_{q^{b}}^{r}\left[a\right]_{q}^{n}\sum_{j_{1},\cdots,j_{r}=0}^{a-1}\left(-1\right)^{{\displaystyle {\textstyle \sum_{l=1}^{r}j_{l}}}}q^{b{\textstyle {\displaystyle {\textstyle \sum_{l=1}^{r}j_{l}}}}}E_{n,q^{a}}^{\left(r\right)}\left(bx+\frac{b}{a}\left(j_{1}+\cdots+j_{r}\right)\right)\\
= & \left[2\right]_{q^{a}}^{r}\left[b\right]_{q}^{n}\sum_{j_{1},\cdots,j_{r}=0}^{b-1}\left(-1\right)^{{\displaystyle {\textstyle \sum_{l=1}^{r}j_{l}}}}q^{a{\displaystyle {\textstyle \sum_{l=1}^{r}j_{l}}}}E_{n,q^{b}}^{\left(r\right)}\left(ax+\frac{a}{b}\left(j_{1}+\cdots+j_{r}\right)\right).
\end{align*}

\end{thm}
$\,$

By (\ref{eq:3}), we easily get
\begin{equation}
E_{n,q}^{\left(r\right)}\left(x+y\right)=\sum_{i=0}^{n}\dbinom{n}{i}q^{xi}E_{i,q}^{\left(r\right)}\left(y\right)\left[x\right]_{q}^{n-i}.\label{eq:9}
\end{equation}

Thus, from (\ref{eq:9}), we have
\begin{align}
 & \sum_{j_{1},\cdots,j_{r}=0}^{a-1}\left(-1\right)^{{\displaystyle {\textstyle \sum_{l=1}^{r}j_{l}}}}q^{b{\displaystyle {\textstyle \sum_{l=1}^{r}j_{l}}}}E_{n,q^{a}}^{\left(r\right)}\left(bx+\frac{b}{a}\left(j_{1}+\cdots+j_{r}\right)\right)\label{eq:10}\\
= & \sum_{j_{1},\cdots,j_{r}=0}^{a-1}\left(-1\right)^{{\displaystyle {\textstyle \sum_{l=1}^{r}j_{l}}}}q^{b{\displaystyle {\textstyle \sum_{l=1}^{r}j_{l}}}}\sum_{i=0}^{n}\dbinom{n}{i}q^{ia\left(\frac{b}{a}{\displaystyle {\textstyle \sum_{l=1}^{r}j_{l}}}\right)}E_{i,q^{a}}^{\left(r\right)}\left(bx\right)\left[\frac{b}{a}\sum_{l=1}^{r}j_{l}\right]_{q^{a}}^{n-i}\nonumber \\
= & \sum_{j_{1},\cdots,j_{r}=0}^{a-1}\left(-1\right)^{{\displaystyle {\textstyle \sum_{l=1}^{r}j_{l}}}}q^{b{\displaystyle {\textstyle \sum_{l=1}^{r}j_{l}}}}\sum_{i=0}^{n}\dbinom{n}{i}q^{\left(n-i\right)b{\textstyle {\displaystyle {\textstyle \sum_{l=1}^{r}j_{l}}}}}E_{n-i,q^{a}}^{\left(r\right)}\left(bx\right)\left[\frac{b}{a}\sum_{l=1}^{r}j_{l}\right]_{q^{a}}^{i}\nonumber \\
= & \sum_{i=0}^{n}\dbinom{n}{i}\left(\frac{\left[b\right]_{q}}{\left[a\right]_{q}}\right)^{i}E_{n-i,q^{a}}^{\left(r\right)}\left(bx\right)\nonumber \\
 & \times\sum_{j_{1},\cdots,j_{r}=0}^{a-1}\left(-1\right)^{{\displaystyle {\textstyle \sum_{l=1}^{r}j_{l}}}}q^{b{\displaystyle {\textstyle \sum_{l=1}^{r}\left(n-i+1\right)j_{l}}}}\left[j_{1}+\cdots+j_{r}\right]_{q^{b}}^{i}\nonumber \\
= & \sum_{i=0}^{n}\dbinom{n}{i}\left(\frac{\left[b\right]_{q}}{\left[a\right]_{q}}\right)^{i}E_{n-i,q^{a}}^{\left(r\right)}\left(bx\right)S_{n,i,q^{b}}^{\left(r\right)}\left(a\right),\nonumber
\end{align}
where
\begin{equation}
S_{n,i,q^{b}}^{\left(r\right)}\left(a\right)=\sum_{j_{1},\cdots,\, j_{r}=0}^{a-1}\left(-1\right)^{\sum_{l=1}^{r}j_{l}}q^{\sum_{l=1}^{r}\left(n-i+1\right)j_{l}}\left[j_{1}+\cdots+j_{r}\right]_{q}^{i}.\label{eq:11}
\end{equation}

From (\ref{eq:10}) and (\ref{eq:11}), we note that
\begin{align}
 & \left[2\right]_{q^{b}}^{r}\left[a\right]_{q}^{n}\sum_{j_{1},\cdots,j_{r}=0}^{a-1}\left(-1\right)^{\sum_{l=1}^{r}j_{l}}q^{b\sum_{l=1}^{r}j_{l}}E_{n,q^{a}}^{\left(r\right)}\left(bx+\frac{b}{a}\left(j_{1}+\cdots+j_{r}\right)\right)\label{eq:12}\\
= & \left[2\right]_{q^{b}}^{r}\sum_{i=0}^{n}\dbinom{n}{i}\left[a\right]_{q}^{n-i}\left[b\right]_{q}^{i}E_{n-i,q^{a}}^{\left(r\right)}\left(bx\right)S_{n,i,q^{b}}^{\left(r\right)}\left(a\right).\nonumber
\end{align}

By the same method as (\ref{eq:12}), we get
\begin{align}
 & \left[2\right]_{q^{a}}^{r}\left[b\right]_{q}^{n}\sum_{j_{1},\cdots,j_{r}=0}^{b-1}\left(-1\right)^{\sum_{l=1}^{r}j_{l}}q^{a\sum_{l=1}^{r}j_{l}}E_{n,q^{b}}^{\left(r\right)}\left(ax+\frac{a}{b}\left(j_{1}+\cdots+j_{r}\right)\right)\label{eq:13}\\
= & \left[2\right]_{q^{a}}^{r}\sum_{i=0}^{n}\dbinom{n}{i}\left[b\right]_{q}^{n-i}\left[a\right]_{q}^{i}E_{n-i,q^{b}}^{\left(r\right)}\left(ax\right)S_{n,i,q^{a}}^{\left(r\right)}\left(b\right).\nonumber
\end{align}

Therefore, by (\ref{eq:12}) and (\ref{eq:13}), we obtain the following
theorem.
\begin{thm}
For $n\ge0$ and $a,\, b\in\mathbb{N}$ with $a\equiv1$ $\textnormal{(mod }2\mathnormal{)}$
and $b\equiv1$ $\textnormal{(mod }2\mathnormal{)}$, we have

\begin{align*}
 & \left[2\right]_{q^{b}}^{r}\sum_{i=0}^{n}\dbinom{n}{i}\left[a\right]_{q}^{n-i}\left[b\right]_{q}^{i}E_{n-i,q^{a}}^{\left(r\right)}\left(bx\right)S_{n,i,q^{b}}^{\left(r\right)}\left(a\right)\\
= & \left[2\right]_{q^{a}}^{r}\sum_{i=0}^{n}\dbinom{n}{i}\left[b\right]_{q}^{n-i}\left[a\right]_{q}^{i}E_{n-i,q^{b}}^{\left(r\right)}\left(ax\right)S_{n,i,q^{a}}^{\left(r\right)}\left(b\right).
\end{align*}

\end{thm}
$\,$

It is not difficult to show that
\begin{align}
 & e^{\left[x\right]_{q}u}\sum_{m_{1},\cdots,m_{r}=0}^{\infty}q^{m_{1}+\cdots+m_{r}}\left(-1\right)^{m_{1}+\cdots+m_{r}}e^{\left[y+m_{1}+\cdots+m_{r}\right]_{q}q^{x}\left(u+v\right)}\label{eq:14}\\
= & e^{-\left[x\right]_{q}v}\sum_{m_{1},\cdots,m_{r}=0}^{\infty}q^{m_{1}+\cdots+m_{r}}\left(-1\right)^{m_{1}+\cdots+m_{r}}e^{\left[x+y+m_{1}+\cdots+m_{r}\right]_{q}\left(u+v\right)}.\nonumber
\end{align}

By (\ref{eq:2}) and (\ref{eq:14}), we get
\begin{align}
 & \sum_{k=0}^{m}\dbinom{m}{k}q^{\left(k+n\right)x}E_{k+n,q}^{\left(r\right)}\left(y\right)\left[x\right]_{q}^{m-k}\label{eq:15}\\
= & \sum_{k=0}^{n}\dbinom{n}{k}E_{m+k,q}^{\left(r\right)}\left(x+y\right)q^{\left(n-k\right)x}\left[-x\right]_{q}^{n-k},\nonumber
\end{align}
where $m,n\ge0$.

Thus, by (\ref{eq:15}), we see that
\begin{align}
 & \sum_{k=0}^{m}\dbinom{m}{k}q^{kx}E_{k+n,q}^{\left(r\right)}\left(y\right)\left[x\right]_{q}^{m-k}\label{eq:16}\\
= & \sum_{k=0}^{n}\dbinom{n}{k}q^{-kx}E_{m+k,q}^{\left(r\right)}\left(x+y\right)\left[-x\right]_{q}^{n-k},\nonumber
\end{align}
where $m,\, n\ge0$.

%%%%%%%%%%%%%%%%%%%%%%%%%%%%%%%%%%%%%%%%%%%%%%%%%%%%%%%%%%%%%%%%%%%%%%%%%%%%%%%%%%%%%%%%%%%%

\bigskip
ACKNOWLEDGEMENTS. This work was supported by the National Research Foundation of Korea(NRF) grant funded by the Korea government(MOE)\\
(No.2012R1A1A2003786 ).
\bigskip

%%%%%%%%%%%%%%%%%%%%%%%%%%%%%%%%%%%%%%%%%%%%%%%%%%%%%%%%%%%%%%%%%%%%%%%%%%%%%%%%%%%%%%%%%%%%%%%%%%%%%%%%%%%%%%%%%%%%%%%%%%%%%%%%%%%%%%%%%%%%%%%%%%%%%%\begin{thebibliography}{99}

\bibliographystyle{amsplain}
\nocite{*}
\bibliography{Symmetric_properties}

$\,$

\noindent \noun{Department of Mathematics, Sogang University, Seoul
121-742, Republic of Korea}

\noindent \emph{E-mail}\noun{ }\emph{address : }\texttt{dskim@sogang.ac.kr}

\noun{$\,$}

\noindent \noun{Department of Mathematics, Kwangwoon University, Seoul
139-701, Republic of Korea}

\noindent \emph{E-mail}\noun{ }\emph{address : }\texttt{tkkim@kw.ac.kr}
\end{document}